\documentclass[11pt,a4paper]{article}

\usepackage{graphicx,latexsym,euscript,makeidx,color,bm}
\usepackage{amsmath,amsfonts,amssymb,amsthm,thmtools,mathtools,mathrsfs,enumerate}
\usepackage[colorlinks,linkcolor=blue,citecolor=red]{hyperref}
\usepackage{tikz}

\usepackage{geometry}
\allowdisplaybreaks[4]
\def\sqr#1#2{{\vcenter{\vbox{\hrule height.#2pt
				\hbox{\vrule width.#2pt height#1pt \kern#1pt \vrule width.#2pt}
				\hrule height.#2pt}}}}
\def\signed #1{{\unskip\nobreak\hfil\penalty50
		\hskip2em\hbox{}\nobreak\hfil#1
		\parfillskip=0pt \finalhyphendemerits=0 \par}}
\def\endpf{\signed {$\sqr69$}}

 \def\sD{\mathscr{D}}  
\def\dbE{\mathbb{E}}   
\def\dbF{\mathbb{F}} \def\sF{\mathscr{F}}  
   
\def\dbH{\mathbb{H}}

 \def\sL{\mathscr{L}}  \def\cL{{\cal L}}
   \def\cM{{\cal M}}

\def\dbP{\mathbb{P}}   
   
\def\dbR{\mathbb{R}}   
\def\dbS{\mathbb{S}}   
   
 \def\sU{\mathscr{U}}

 \def\sX{\mathscr{X}}

\def\ss{\smallskip}             \def\hb{\hbox}
\def\ms{\medskip}              
          \def\lan{\langle}       
\def\h{\widehat}     \def\ran{\rangle}       \def\tr{\hb{tr$\,$}}
\def\q{\quad}               
\def\qq{\qquad}             
\def\no{\noindent}         \def\sc{\scriptscriptstyle}
         \def\scT{\sc T}
    \def\Blan{\Big\lan\!\!} \def\ds{\displaystyle}
\def\rf{\eqref}       \def\Bran{\!\!\Big\ran} 
\def\cd{\cdot}        \def\({\Big(}           
       \def\){\Big)}           \def\les{\leqslant}
\def\wt{\widetilde}   \def\[{\Big[}           \def\ges{\geqslant}
  \def\]{\Big]}           
\def\a{\alpha}       \def\l{\lambda}       
\def\b{\beta}        \def\t{\tau}            
\def\d{\delta}           \def\G{\varGamma}   
\def\e{\varepsilon}      \def\L{\varLambda}
\def\f{\varphi}           
       \def\si{\sigma}    \def\Si{\varSigma}
\def\i{\infty}             \def\Th{\varTheta}
\def\be{\begin{equation}}
\def\bel{\begin{equation}\label}
\def\ee{\end{equation}}
\def\bea{\begin{eqnarray}}
\def\eea{\end{eqnarray}}
\def\bt{\begin{theorem}}
\def\et{\end{theorem}}
\def\bc{\begin{corollary}}
\def\ec{\end{corollary}}
\def\bl{\begin{lemma}}
\def\el{\end{lemma}}
\def\bp{\begin{proposition}}
\def\ep{\end{proposition}}
\def\br{\begin{remark}}
\def\er{\end{remark}}
\def\ba{\begin{array}}
\def\ea{\end{array}}
\def\bd{\begin{definition}}
\def\ed{\end{definition}}
\def\5n{\negthinspace \negthinspace \negthinspace \negthinspace \negthinspace }
\def\4n{\negthinspace \negthinspace \negthinspace \negthinspace }
\def\3n{\negthinspace \negthinspace \negthinspace }
\def\2n{\negthinspace \negthinspace }
\def\1n{\negthinspace }

\newtheoremstyle{thry}
{}      
{}      
{\sl}   
{}      
{\bf}   
{.}     
{.5em}  
{}      
\theoremstyle{thry}

\newtheorem{theorem}{Theorem}[section]
\newtheorem{proposition}[theorem]{Proposition}
\newtheorem{corollary}[theorem]{Corollary}
\newtheorem{lemma}[theorem]{Lemma}

\theoremstyle{definition}
\newtheorem{definition}[theorem]{Definition}
\newtheorem{example}[theorem]{Example}

\theoremstyle{remark}
\newtheorem{remark}[theorem]{Remark}

\def\punct{}
\newtheoremstyle{dotless}{}{}{\rm}{}{\bf}{\punct}{.5em}{}
\theoremstyle{dotless}

\makeatletter
   
   \@addtoreset{equation}{section}
   \newcommand{\setword}[2]{%
   \phantomsection
   #1\def\@currentlabel{\unexpanded{#1}}\label{#2}%
   }
\makeatother

\def\3n{\negthinspace \negthinspace \negthinspace }
\def\2n{\negthinspace \negthinspace }
\def\1n{\negthinspace }
\def\bel{\begin{equation}\label}

\def\dbE{\mathbb{E}}
\def\dbF{\mathbb{F}}

\def\dbH{\mathbb{H}}

\def\dbP{\mathbb{P}}

\def\dbR{\mathbb{R}}
\def\dbS{\mathbb{S}}

\def\sD{\mathscr{D}}

\def\sF{\mathscr{F}}

\def\sL{\mathscr{L}}

\def\sU{\mathscr{U}}

\def\sX{\mathscr{X}}

\def\ds{\displaystyle}

\def\ns{\noalign{\ss}}
\def\a{\alpha}
\def\b{\beta }

\def\d{\delta}
\def\e{\varepsilon}

\def\l{\lambda}

\def\si{\sigma}
\def\t{\tau}
\def\f{\varphi}

\def\i{\infty}
\def\G{\Gamma}

\def\Th{\Theta}
\def\L{\Lambda}
\def\Si{\Sigma}

\def\O{\Omega}
\def\cL{{\cal L}}
\def\cM{{\cal M}}

\def\BS{{\bf S}}

\def\G{\Gamma}

\def\Th{\Theta}
\def\L{\Lambda}
\def\Si{\Sigma}

\def\O{\Omega}

\def\no{\noindent}

\def\ss{\smallskip}
\def\ms{\medskip}

\def\q{\quad}
\def\qq{\qquad}
\def\hb{\hbox}

\def\h1{\outline{$1$}}

\def\hSi{\outline{$\Sigma$}}

\def\hTh{\outline{$\Theta$}}

\def\hh1{\outline{$1$}}
\def\hh2{\outline{$2$}}
\def\hh3{\outline{$3$}}
\def\hh4{\outline{$4$}}
\def\hh5{\outline{$5$}}
\def\hh6{\outline{$6$}}
\def\hh7{\outline{$7$}}
\def\hh8{\outline{$8$}}
\def\hh9{\outline{$9$}}
\def\hh0{\outline{$0$}}

\def\limsup{\mathop{\overline{\rm lim}}}

\def\lan{{\langle}}
\def\ran{{\rangle}}

\def\pa{\partial}
\def\h{\widehat}
\def\wt{\widetilde}

\def\cd{\cdot}
\def\cds{\cdots}

\def\tr{\hbox{\rm tr$\,$}}

\def\les{\leqslant}
\def\ges{\geqslant}

\def\({\Big (}
\def\){\Big )}
\def\[{\Big[}
\def\]{\Big]}
\def\lan{\langle}
\def\ran{\rangle}

\def\bde{\begin{definition}\label}
	\def\ede{\end{definition}}
\def\bel{\begin{equation}\label}
		\def\ee{\end{equation}}
	\def\bt{\begin{theorem}\label}
		\def\et{\end{theorem}}
	\def\bc{\begin{corollary}\label}
		\def\ec{\end{corollary}}
	\def\bl{\begin{lemma}\label}
		\def\el{\end{lemma}}
	\def\bp{\begin{proposition}\label}
		\def\ep{\end{proposition}}
	\def\bex{\begin{example}\label}
		\def\ex{\end{example}}
	\def\bas{\begin{assumption}}
		\def\eas{\end{assumption}}
	\def\br{\begin{remark}\label}
		\def\er{\end{remark}}
	\def\ba{\begin{array}}
		\def\ea{\end{array}}
	\def\ed{\end{document}}

\def\rf{\eqref}

\def\square#1{\vbox{\hrule\hbox{\vrule height#1%
			\kern#1\vrule}\hrule}}
\def\rectangle#1#2{\vbox{\hrule\hbox{\vrule height#1%
			\kern#2\vrule}\hrule}}

\font\tenbb=msbm10 \font\sevenbb=msbm7 \font\fivebb=msbm5

\newfam\bbfam
\scriptscriptfont\bbfam=\fivebb \textfont\bbfam=\tenbb
\scriptfont\bbfam=\sevenbb

\makeatletter

\@addtoreset{equation}{section}
\makeatother

\usepackage{xcolor}
\makeatletter

\input pdf-trans
\newbox\qbox
\def\usecolor#1{\csname\string\color@#1\endcsname\space}
\newcommand\bordercolor[1]{\colsplit{1}{#1}}
\newcommand\fillcolor[1]{\colsplit{0}{#1}}
\newcommand\outline[1]{\leavevmode%
	\def\maltext{#1}%
	\setbox\qbox=\hbox{\maltext}%
	\boxgs{Q q 2 Tr \thickness\space w \fillcol\space \bordercol\space}{}%
	\copy\qbox%
}
\makeatother
\newcommand\colsplit[2]{\colorlet{tmpcolor}{#2}\edef\tmp{\usecolor{tmpcolor}}%
	\def\tmpB{}\expandafter\colsplithelp\tmp\relax%
	\ifnum0=#1\relax\edef\fillcol{\tmpB}\else\edef\bordercol{\tmpC}\fi}
\def\colsplithelp#1#2 #3\relax{%
	\edef\tmpB{\tmpB#1#2 }%
	\ifnum `#1>`9\relax\def\tmpC{#3}\else\colsplithelp#3\relax\fi
}
\bordercolor{black}
\fillcolor{white}
\def\thickness{.3}

\def\hSi{\outline{$\Sigma$}}

\def\hTh{\outline{$\Theta$}}

\def\h1{\outline{$1$}}

\begin{document}

\title{\bf Turnpike Property of Stochastic Linear-Quadratic Optimal Control Problems in Large Horizons with Regime Switching I: Homogeneous Cases}
		
		\author{Hongwei Mei\footnote{ Department of Mathematics and Statistics, Texas Tech University, Lubbock, TX 79409, USA; email: {\tt hongwei.mei@ttu.edu}. This author is partially supported by Simons Travel Grant MP-TSM-00002835.},~~~Rui Wang\footnote{ Department of Mathematics and Statistics, Texas Tech University, Lubbock, TX 79409, USA; email: {\tt rui-math.wang@ttu.edu}},~~~\text{and}~~~
			Jiongmin Yong\footnote{Department of Mathematics, University of Central Florida, Orlando, FL 32816, USA; email: {\tt jiongmin.yong@ucf.edu}. This author is supported in part by NSF Grant DMS-2305475.}  }
				
\maketitle

\no{\bf Abstract.}
This paper is concerned with optimal control problems for a linear homogeneous stochastic differential equation having regime switching with purely quadratic functional in the large time horizons. We establish the so-called turnpike properties for the optimal pairs. The key is to prove a proper convergence of the solutions to the differential Riccati equations to the algebraic Riccati equation.  Even for the problems without regime switchings, our result provides a refined estimate compared to those in the previous literature, which also provides a new tool for further research.

\ms
\no{\bf Keywords.}
Turnpike property, regime switching, stochastic optimal control, linear-quadratic, static optimization,
stabilizability, Riccati equation.

\ms
\no{\bf AMS 2020 Mathematics Subject Classification.}  49N10, 93E15, 93E20.

\section{Introduction}\label{Sec:Intro}

Let $(\Omega,\sF,\dbP)$ be a complete probability space on which a standard one-dimensional Brownian motion $W=\{W(t);\,t\ges 0\}$ and a  Markov chain $\a(\cd)$ with a finite state space $\cM=\{1,2,3,\cds,m_0\}$ are defined, for which they are assumed to be independent. The generator of $\a(\cd)$ is denoted by $(\l_{\imath\jmath})_{m_0\times m_0}$ (see below for details). We now denote by $\dbF^W=\{\sF^W_{t}\}_{t\ges0}$ (resp. $\dbF^\a=\{\sF_{t}^\a\}_{t\ges0}$, $\dbF=\{\sF_{t}\}_{t\ges0})$ the usual augmentation of the natural filtration generated by $W(\cd)$ (resp. by $\a(\cd)$, and by $(W(\cd),\a(\cd))$).
Consider the following {\it state equation} which is a controlled linear homogeneous stochastic differential equation (SDE, for short), with regime switchings:
\bel{state}\left\{\2n\ba{ll}
\ns\ds dX(s)=\big\{A(\a(s))X(s)+B(\a(s))u(s)\big\}ds\\
\ns\ds\qq\qq+\big\{C(\a(s))X(s)+D(\a(s))u(s)\big\}dW(s),\q s\in[t,T)],\\
\ns\ds X(t)=x,\q\a(t)=\imath,\ea\right.\ee
where $[t,T)]$ stands for $[t,T]$ if $T<\i$ and $[t,\i)$ if $T=\i$. We adopt the following basic assumption for the {\it coefficients} of the state equation \rf{state}:

\ms

{\bf(A1)} Let $A,C:\cM\to\dbR^{n\times n}$ and $B,D:\cM\to\dbR^{m\times n}$ be measurable.

\ms

Since $\cM$ is finite, we automatically have the boundedeness of these coefficients. In the case that $\cM$ is a singleton, our following problems will be all reduced to those without regime switchings. Now, we introduce the spaces (with $0<T\les\i$)
$$\ba{ll}
\ns\ds L^2_\dbF(t,T;\dbH)=\Big\{\f:[t,T)\times\O\to\dbH\bigm|\f(\cd)\hb{ is $\dbF$-progressively measurable, }\\
\ns\ds\qq\qq\qq\qq\qq\qq\qq\qq\qq\qq\dbE\int_t^T|\f(s)|_{\dbH}^2ds<\i\Big\},\\
\ns\ds L^2_{\dbF^\a}(t,T;\dbH)=\Big\{\f:[0,T)\times\O\to\dbH\bigm|\f(\cd)\hb{ is $\dbF^\a$-progressively measurable, }\\
\ns\ds\qq\qq\qq\qq\qq\qq\qq\qq\qq\qq\dbE\int_t^T|\f(s)|_{\dbH}^2ds<\i\Big\},\ea$$

\ms

In \rf{state}, $(t,x,\imath)\in[0,\i)\times\dbR^n\times\cM\equiv\sD$ is called the {\it initial triple}, $X(\cd)$ is called the {\it state process}, and $u(\cd)$, called a {\it control}, is selected from the space
$$\sU[t,T]=L_\dbF^2(t,T;\dbR^m).$$
It is clear that for each $0\les t<T\les\i$, any $(t,x,\imath)\in\sD$ and $u(\cd)\in\sU[t,T]$, under (A1), state equation \rf{state} admits a unique solution $X(\cd)=X(\cd\,;t,x,\imath;u(\cd))\in\sX[t,T]$ on $[t,T)]$, with the set of {\it feasible state processes}:
$$\ba{ll}
\ns\ds\sX[t,T]\equiv\Big\{X:[t,T)]\to\dbR^n\bigm|X(\cd)\hb{ is $\dbF$-adaped, has continuous paths, and}\\
\ns\ds\qq\qq\qq\qq\qq\qq\dbE\[\sup_{s\in[t,S]}|X(s)|^2\]<\i,\q \forall s<T\Big\}.\ea$$
Now, for $0\les t<T\les\i$, we let the set of {\it admissible state processes} be
$$\sX_{ad}[t,T]=\sX[t,T]\cap L^2_\dbF(t,T;\dbR^n),$$
and the set of {\it admissible controls} be
$$\sU_{ad}^{t,x,\imath}[t,T]=\Big\{u(\cd)\in\sU[t,T]\bigm|
X(\cd\,;t,x,\imath;u(\cd))\in\sX_{ad}[t,T]\Big\}.$$
Clearly,
$$\sU_{ad}^{t,x,\imath}[t,T]=\sU[t,T],\qq0\les t<T<\i,\q(t,x,\imath)\in\sD.$$
But, in general,
\bel{U ne U}\sU_{ad}^{t,x,\imath}[t,\i]\subsetneq\sU[t,\i],\ee
Namely, for some $u(\cd)\in\sU[t,\i]$ (and some $(t,x,\imath)\in\sD$), we might have $X(\cd\,;t,x,\imath;u(\cd))\notin\sX_{ad}[t,\i]$ in general. Rigorously speaking, the state-control pair $(X(\cd),u(\cd))$ depends on the Markov chain $\a(\cd)$, since the former are merely $\dbF$-progressively measurable.

\ms

To measure the performance of an admissible control $u(\cd)\in\sU_{ad}^{t,x,\imath}[t,T]$, we introduce the following cost functional
\bel{cost[t,T]}J_{\scT}(t,x,\imath;u(\cd))=\dbE\(\int_t^Tg
(s,X(s),\a(s),u(s))ds\),\ee
where
\bel{g}g(s,x,\imath,u)={1\over2}\Blan\begin{pmatrix}Q(\imath)& S(\imath)^\top \\ S(\imath)&R(\imath)\end{pmatrix}\begin{pmatrix}x\\ u\end{pmatrix},\begin{pmatrix}x\\ u\end{pmatrix}\Bran.\ee
Here, the superscript $\top$ denotes the transpose of matrices; $\lan\cd\,,\cd\ran$ denotes the inner product of two vectors in proper spaces which can be identified from the context. Next, we denote $\dbS^n$, $\dbS^n_+$ and $\dbS^n_{++}$ to be the sets of all $(n\times n)$ symmetric, positive semi-definite, and positive definite matrices, respectively. Now it is natural to consider the following optimal control problem for $0\les t<T\les\i$.

\ms

{\bf Problem (LQ)$_{\sc t,\scT}$.} For a given initial triple $(t,x,\imath)\in\sD$, find a control $\bar u_{\scT}^{t,x,\imath}(\cd)\in\sU_{ad}^{t,x,\imath}[t,T]$ such that
\bel{opt-u3}J_{\scT}(t,x,\imath;\bar u_{\scT}^{t,x,\imath}(\cd))=\inf_{u(\cd)\in\sU_{ad}^{t,x,\imath}
[t,T]}J_{\scT}(t,x,\imath;u(\cd))\equiv V_{\scT}(t,x,\imath).\ee

\ms

The above problem is referred to as a (homogeneous) {\it linear-quadratic} (LQ, for short) {\it optimal control problem} over the finite or infinite time horizon with regime switchings (see \cite{Zhang-2021, Mei-Wei-Yong-2025} for examples). We call $\bar u^{t,x,\imath}_{\scT}(\cd)$ an {\it optimal control process}, $\bar X^{t,x,\imath}_{\scT}(\cd)$ the corresponding {\it optimal state process}, and $(\bar X^{t,x,\imath}_{\scT}(\cd),\bar u^{t,x,\imath}_{\scT}(\cd))$ the {\it optimal pair} of Problem (LQ)$_{\sc t,\scT}$, respectively. In addition, we call $V_T(t,x,\imath)$ the {\it optimal value} of Problem (LQ)$_{\sc t,\scT}$, and thus as $(t,x,\imath)\in\sD$ varies, it is called the {\it value function} of the problem.
Under some general mild assumptions, it can be proved that Problem (LQ)$_{\sc t,\scT}$ admits a unique optimal control $\bar u_{\scT}^{t,x,\imath}(\cd)\in\sU_{ad}^{t,x,\imath}[t,T]$ with the optimal trajectory $\bar X_{\scT}^{t,x,\imath}(\cd)\equiv X(\cd\,;t,x,\imath;\bar u_{\scT}^{t,x,\imath}(\cd))\in L^2_\dbF(t,T;\dbR^n)$. For the cases without regime switchings, under proper conditions (see the next section for similar conditions), people found that for some stochastic processes $(\bar X_\i(\cd),\bar u_\i(\cd))$, and some constants $\b,K>0$, all are independent of $0<T<\i$, (in what follows, $K>0$ will be a generic constant which can be different from line to line) such that
\bel{STP0}\dbE(|\bar X_{\scT}^{0,x}(s)-\bar X_\i(s)|^2+|\bar u_{\scT}^{0,x}(s)-\bar u_\i(s)|^2)\les K(e^{-\b(T-s)}+e^{-\b s}),\q\forall t\in[0,T],\ee
where $(\bar X_{\scT}^{0,x}(\cd),\bar u_{\scT}^{0,x}(\cd))$ is the optimal pair of the corresponding Problem (LQ)$_{\sc0,\scT}$ on time horizon $[0,T]$ (with $0<T<\i$, without regime switchings). Such an asymptotic behavior of the optimal pair $(\bar X_{\scT}^{0,x}(\cd),\bar u_{\scT}^{0,x}(\cd))$ as $T\to\i$ is called the {\it strong turnpike property} (STP, for short) of Problem (LQ)$_{\sc0,\scT}$ (without regime switchings). A key interesting
feature of \rf{STP0} is that, for a very large $T$, in the middle part of interval $[0,T]$, the optimal pair $(\bar X^{0,x}_{\scT}(\cd),\bar u^{0,x}_{\scT}(\cd))$ is very close to a given $T$-independent pair $(\bar X^{0,x}_{\sc\i}(\cd),\bar u^{0,x}_{\sc\i}(\cd))$. Thus, the latter is a good approximation of the former in the middle part of $[0,T]$.

\ms

Turnpike phenomenon was realized by Ramsey in 1928 (\cite{Ramsey-1928}). Formal investigations on turnpike property (for deterministic economics systems) can be dated back to von Neumann \cite{Neumann-1945} in 1945. Later, in 1958, Dorfman--Samuelson--Solow (\cite{Dorfman-Samuelson-Solow-1958}) coined the name {\it turnpike property}, intuitively suggested by the highway system of the United States. Since then the turnpike property has been found to hold for a large class of (deterministic, finite or infinite dimensional) optimal control problems. Numerous relevant results can be found in \cite{McKenzie-1976,Carlson-Haurie-Leizarowitz-1991,Damm-Grune-Stieler-Worthmann-2014,
Trelat-Zuazua-2015,Zuazua-2017,Grune-Guglielmi-2018,Zaslavski-2019,
Lou-Wang-2019,Breiten-Pfeiffer-2020,Sakamoto-Zuazua-2021,Faulwasser-Grune-2022} and the references cited therein. Turnpike properties for stochastic optimal control problems were not touched for a long time. In the book \cite{Carlson-Haurie-Leizarowitz-1991}, some random ODEs (i.e., no It\^o's integral was involved) with certain random jumps, whose rate was governed by a Markov type chain, the corresponding turnpike property was studied. At about the same time, certain stability for a finite time horizon multi-person discrete stochastic game was investigated and using the idea of turnpike property, it was shown the existence of an equilibrium for a stationary (discrete random) games (\cite{Marimon-1989}). A systematic investigation for continuous time stochastic optimal LQ control problems was begun by the work of Sun--Wang--Yong in the early of 2020 (\cite{Sun-Wang-Yong-2022}), followed by the works  \cite{Conforti-2023,Chen-Luo-2023,Sun-Yong-2024,Sun-Yong-2024b,
Jian-Jin-Song-Yong-2024,Schiessl-Baumann-Faulwasser-Grune-2024,
Bayraktar-Jian-2025}. Stochastic control systems with regime switchings have broad applications in various areas. Then it is naturally to ask do we have turnpike type properties for stochastic optimal control with system containing regime switchings? The main purpose of this paper is to give a positive answer and refine the known results. More precisely, we will establish the following (compare with \rf{STP0}): there exist constants $\d,K>0$ and some positive function $h(\cd)$ on $ [0,\i)$, independent of $0<T<\i$ such that 
\bel{new}\ba{ll}
\ns\ds\dbE\Big(|\bar X_{\scT}^{t,\imath,x_{\scT}}(s)-\bar X_{\sc\i}^{t,\imath,x_\i}(s)|^2+|\bar u_{\scT}^{t,\imath,x_T}(s)-\bar u_{\sc\i}^{t,\imath,x_\i}(s)|^2\Big)\\
\ns\ds\les Ke^{-\d(s-t)}|x_{\scT}-x_{\sc\i}|^2+K
e^{-2\d(T-s)}\Big(e^{-\d(s-t)}|x_{\scT}|^2+h(s-t)\Big),\qq s\in[t,T],\ea\ee
with $(t,x_{\scT},\imath),~(t,x_{\sc\i},\imath)\in\sD$ being two possibly different initial triples for the optimal pairs $(\bar X_{\scT}^{t,\imath,x_{\scT}}(\cd),\bar u_{\scT}^{t,\imath,x_T}(\cd))$ and $(\bar X_{\sc\i}^{t,\imath,x_\i}(\cd),\bar u_{\sc\i}^{t,\imath,x_\i}(\cd))$ of Problems (LQ)$_{\sc t,T}$ and (LQ)$_{\sc t,\i}$, respectively. In the homogeneous cases studied in the paper, it turns out that $h(\cd)=0.$ Observing from this,  if we take $x_{\scT}=x_{\sc\i}=x$, $t=0$ and use the well-known fact $2ab\les a^2+b^2$, by some possibly different constants, one can easily recover the results from \rf{STP0}.

\ms

In this paper, we will focus on the results for homogeneous cases only and those for non-homogeneous cases will involve other techniques and will be presented in a subsequent paper (with a different function $h(\cd)$).  The key of this paper is to establish a proper convergence of the solutions to certain differential Riccati equations to an algebraic Riccati equation. This will also play an essential role for non-homogeneous cases (see the subsequent paper). The rest of the paper is arranged as follows. Section \ref{sec:sta} recalls some classical results on linear control system with regime switchings, including stabilization and the extended It\^o's formula for our system with regime switchings. Section \ref{sec:cost} is devoted to the characterization of the cost functional and the value function. Then the turnpike property for homogeneous cases is established in Section \ref{sec:hom}, where we find a more straightforward proof for the exponential convergence of the solutions to the differential Riccati equations different from that in \cite{Sun-Wang-Yong-2022}. Finally, some concluding remarks are made in Section \ref{sec:con}.

\section{A General Consideration}\label{sec:sta}

\subsection{Stabilization and dissipating strategy}

First, we recall the generator $(\l_{\imath\jmath})_{m_0\times m_0}\in\dbR^{m_0\times m_0}$ of Markov chain $\a(\cd)$ with a finite state space $\cM$, which is a real matrix so that the following hold:
\bel{q-prop}\l_{\imath\jmath}>0,\q\imath\ne\jmath,\qq\sum_{\jmath=1}^{m_0}
\l_{\imath\jmath}=0,\q\imath\in\cM.\ee
For any measurable $\Si:\cM\to\dbS$, define
\bel{L}\L[\Si](\imath)=\sum_{j\in\cM}\l_{\imath\jmath}\Si(\jmath).\ee
Let
$$\ba{ll}
\ns\ds\hTh=\Big\{\Th:\cM\to\dbR^{m\times n}\bigm|\Th(\cd)\hb{ is measurable}\Big\},\\
\ns\ds\hSi=\Big\{\Si:\cM\to\dbS^n_{++}\bigm|\Si(\cd)\hb{ is measurable}\Big\}.\ea$$
Now we first consider the following linear SDE with a regime switching governed by a Markov chain:
\bel{state3}\left\{\2n\ba{ll}
\ns\ds dX(s)=A(\a(s))X(s)ds+C(\a(s))X(s)dW(s),\qq s\in[t,\i),\\
\ns\ds X(t)=x,\q\a(t)=\imath.\ea\right.\ee
The above system is denoted by $[A,C]$. Under (A1), such a system is well-posed. If $X(\cd)\equiv X(\cd\,;t,x,\imath)$ is the solution of the above corresponding to $(t,x,\imath)\in\sD$. Then for any $\Si(\cd)\in\hTh$, it holds
\bel{SiXX}d\lan\Si(\a(s))X(s),X(s)\ran=\Big\lan\(\L[\Si]+\Si A+A^\top\Si+C^\top\Si C\)(\a(s))X(s),X(s)\Big\ran.\ee
This is an application of the extended It\^o's formula, see the next subsection for details.

\bde{stable} \rm (i) System $[A,C]$ is said to be {\it stable} if for any $(t,x,\imath)\in\sD$, $X(\cd\,;t,x,\imath)\in\sX_{ad}[t,\i]$.

\ms

(ii) System $[A,C]$ is said to be {\it dissipative} if one could find a $\Si(\cd)\in\hSi$ and a $\d>0$ so that
\bel{dissipative}\(\L[\Si]+\Si A+A^\top\Si+C^\top\Si C\)(\imath)\les-\d\Si(\imath),\qq\imath\in\cM.\ee

\ede

\ms

The constant $\d>0$ in \rf{dissipative} is intrinsically determined by (the stability of) the linear SDE \rf{state3}, namely, this $\d$ cannot be arbitrarily large. For example, when $d=1$ and $\cM=\{1\}$, then \rf{dissipative} reads
$$(2A+C^2)\Si\les-\d\Si.$$
Hence, $0<\d\les-2A-C^2$ (for which $2A+C^2$ is necessarily negative). Now, we consider state equation \rf{state}, which is denoted by $[A,C;B,D]$. Next, for any $\Th\in \hTh$, we write
\bel{A^Th}A^\Th(\imath):=A(\imath)+B(\imath)\Th(\imath),\q C^\Th(\imath):=C(\imath)+D(\imath)\Th(\imath).\ee
The following definition is adopted from \cite{Mei-Wei-Yong-2025}.

\bde{Def-stab-1} (i) System $[A,C;B,D]$ is said to be {\it stabilizable} if one can find a map $\Th(\cd)\in\hTh$, so that for $[A^\Th,C^\Th]$ is stable. In this case, the map $\Th(\cd)$ is called a {\it stabilizer} of $[A,C;B,D]$. The set of all possible stabilizers of system $[A,C;B,D]$ is denoted by $\BS[A,C; B,D]$.

\ms

(ii) The map $\Th(\cd)\in\hTh$ is called a {\it dissipating strategy} of system $[A,C;B,D]$ if there exists a $\d>0$ and a $\Si(\cd)\in\hSi$ such that \rf{dissipative} holds with $[A,C]$ replaced by $[A^\Th,C^\Th]$.

\ede

Since the state space of Markov chain is finite, the following is true (see \cite{Mei-Wei-Yong-2025}, Proposition 3.7).

\bp{} System $[A,C;B,D]$ is stabilizable if and only if it admits a dissipating strategy.

\ep

It is known that stabilizability of $[A,C;B,D]$ is necessary for studying LQ problems in an infinite time horizon even for the LQ problems without regime switchings (\cite{Sun-Yong-2020, Mei-Wei-Yong-2025}). Thus, we introduce the following assumption.

\ms

{\bf (A2)} System $[A,C;B,D]$ is stabilizable, i.e., $\BS[A,C;B,D]\ne\varnothing$.

\ms

Now, let $\Th^*(\cd)\in\BS[A,C;B,D]$ be fixed, $v(\cd)\in\sU[0,T]$ and consider the following controlled system:
\bel{state-closed}\left\{\2n\ba{ll}
\ns\ds dX(s)=\big\{A^{\Th^*}(\a(s))X(s)+B(\a(s))v(s)\big\}ds\\
\ns\ds\qq\qq+\big\{C^{\Th*}(\a(s))X(s)+D(\a(s))v(s))\big\}
dW(s),\q s\in[t,T)],\\
\ns\ds X(t)=x,\q\a(t)=\imath,\ea\right.\ee
which is nothing but the state equation \rf{state} under control
\bel{feedback}u(s)=\Th^*(\a(s))X(s)+v(s),\qq s\in[t,T)].\ee
We denote the unique solution of \rf{state-closed} by $X(\cd)=X(\cd\,;t,x,\imath,(\Th^*(\cd),v(\cd)))$.
The following result is Proposition 4.5 in \cite{Mei-Wei-Yong-2025}.

\bp{equivalence} \sl Let {\rm(A1)} hold. Let $\Th^*(\cd)\in\BS[A,C;B,D]$. Then
\bel{su}\sU_{ad}^{t,x,\imath}[t,\i]=\Big\{\Th^*(\a(\cd))X(\cd\,;
t,x,\imath,(\Th^*(\cd),v(\cd)))+v(\cd)\bigm|v(\cd)\in\sU[t,\i]\Big\}.\ee

\ep

We note that the right-hand side in \eqref{su} is directly independent of the choice of $(t,x,\imath)\in\sD$, which has a clear advantage. Now, fix $\Th^*(\cd)\in\BS[A,C; B,D]$, which means 
\bel{stable-1}0\in\BS[A^{\Th^*},C^{\Th^*};B,D],\ee
or $[A^{\Th^*},C^{\Th^*}]$ is stable. By taking the control as in \rf{feedback}, we have
\bel{costTh}\ba{ll}
\ns\ds J_{\sc\i}(t,x,\imath;u(\cd))=J_{\sc\i}(t,x,\imath;\Th^*(\a(\cd))X(\cd)+v(\cd))\\
\ns\ds\q=\dbE\int_t^\i\[\lan Q(\a(s))X(s),X(s)\ran+2\lan S(\a(s))X(s),\Th^*(\a(s))X(s)+v(s)\ran\\
\ns\ds\qq\qq\qq+\lan R(\a(s))\Th^*(\a(s))X(s)+v(s),\Th^*(\a(s))X(s)+v(s)\ran\]ds\\
\ns\ds\q=\dbE\int_t^\i\[\lan Q^{\Th^*}(\a(s))X(s),X(s)\ran+2\lan S^{\Th^*}(\a(s))X(s),v(s)\ran\\
\ns\ds\qq\qq\qq+\lan R(\a(s))v(s),v(s)\ran\]ds\\
\ns\ds\q\equiv J^*_{\sc\i}(t,x,\imath;v(\cd)),\ea\ee
where
$$\ba{ll}
\ns\ds Q^{\Th^*}(\a(s))\1n=\1n Q(\a(s))\1n+\1n\Th^*(\a(s))^\top\1n S(\a(s))\1n+\1n S(\a(s))^\top\Th^*(\a(s))\1n+\1n\Th^*(\a(s))^\top\1n R(\a(s))\Th^*(\a(s)),\\
\ns\ds S^{\Th^*}(\a(s))=S(\a(s))+R(\a(s))\Th^*(\a(s)).\ea$$
Then under (A1) and (A2), for any fixed $\Th^*(\cd)\in\BS[A,C;B,D]$, we can state a new LQ problem as follows.

\ms

{\bf Problem (LQ)$_{\sc t,\sc\i}^*$.}  For any $(t,x,\imath)\in\sD$, find a $\bar v_{\sc\i}^{t,x,\imath}(\cd)\in\sU[t,T]$ such that
$$V^*_{\sc\i}(t,x,\imath)=J^*_{\sc\i}(t,x,\imath;\bar v_{\sc\i}^{t,x,\imath}(\cd))=\inf_{v(\cd)\in\sU[t,\i]}J^*_{\sc\i}
(t,x,\imath;v(\cd)).$$

\ms

It is clear that the optimal pair of the above Problem (LQ)$_{\sc t,\i}^*$ (if it exists) is given by $(\bar X^{t,x,\imath}_{\sc\i}(\cd),\bar v_{\sc\i}^{t,x,\imath}(\cd))$, where
\bel{v=u}\bar v_{\sc\i}^{t,x,\imath}(\cd)=\bar u_{\sc\i}^{t,x,\imath}(\cd)-\Th^*(\a(\cd))\bar X_{\sc\i}^{t,x,\imath}(\cd).\ee
with $(\bar X^{t,x,\imath}_{\sc\i}(\cd),\bar u_{\sc\i}^{t,x,\imath}(\cd))$ being the optimal pair of Problem ({LQ)$_{\sc t,\i}$.

\ms

From above, we have seen two facts:

\ms

$\bullet$ Without the loss of generality, instead of (A2), we may assume the following in the sequel:

\ms

{\bf(A2)$'$} Let $0\in\BS[A,C;B,D]$, or $[A,C]$ is stable.

\ms

\no Otherwise we can work with $J_{\sc\i}^*(t,x,\imath;v(\cd))$ instead of $J_{\sc\i}(t,x,\imath;u(\cd))$.

\ms

$\bullet$ For any given $\Th^*(\cd)\in\BS[A,C;B,D]$, Problems (LQ)$_{\sc t,\i}$ and (LQ)$_{\sc t,\i}^*$ are equivalent. Here, we should note that although the cost functionals in these two problems look quite different, under the relation \rf{v=u}, we have the equality \rf{costTh}.

\ms

Finally, the same argument can be applied to Problems (LQ)$_{\sc t,T}$ and formulate the corresponding (LQ)$_{\sc t,T}^*$ with $T<\i$. Then we have equivalence between these two problems. Note that in the case $T<\i$, (A2) or (A2)$'$ is not needed.

\subsection{Extended It\^o's formula}

Let us recall the extended It\^o formula for SDEs with regime switchings (quoted from \cite{Skorokhod-1989}, p.104; and \cite{Mao-1999}, pp.47--48): Consider a general SDE
\bel{SDE}dX(t)=b(t,X(t),\a(t))dt+\si(t,X(t),\a(t))dW(t),\ee
with a Markov chain $\a(\cd)$ independent of the Brownian motion $W(\cd)$, whose generator is $\L=(\l_{\imath\jmath})_{m_0\times m_0}$.

\bp{} \sl Let all the required conditions hold. Let $(X(\cd),\a(\cd))$ be a solution to \rf{SDE}. Then for any stopping times $0\les\t_1\les\t_2<\i$,
\bel{Ito}\dbE f(\t_2,X(\t_2),\a(\t_2))-\dbE f(\t_1,X(\t_1),\a(\t_1))=\dbE\int_{\t_1}^{\t_2}\sL f(s,X(s),\a(s))ds,\ee
where
$$\ba{ll}
\ns\ds\sL f(t,x,\imath)=f_t(t,x,\imath)+f_x(t,x,\imath)b(x,\imath)+\tr[f_{xx}(t,x,
\imath)\si(t,x,\imath)
\si(t,x,\imath)]\\
\ns\ds\qq\qq\qq+\L(x)[f(t,x,\cd)](\imath)\\
\ns\ds\qq\equiv{\pa f(t,x,\imath)\over\pa t}+\sum_{i=1}^nb_i(t,x,\imath){\pa f(t,x,\imath)\over\pa x_i}+{1\over2}\sum_{i,j=1}^n\sum_{\ell=1}^d\si_{id}(t,x,\imath)
\si_{jd}(t,x,\imath){\pa^2f(t,x,\imath)\over\pa x_i\pa x_j}\\
\ns\ds\qq\q+\sum_{\jmath\in\cM}\l_{\imath\jmath}[f(t,x,\jmath)
-f(t,x,\imath)].\ea$$

\ep

Now, if we take
$$f(t,x,\imath)={1\over2}\lan P(t,\imath)x,x\ran,$$
for a suitable $P(\cd\,,\cd)$, then the following corollary holds.

\bc{cor2.7} \sl Let {\rm(A1)} hold. Let $(X(\cd),u(\cd))$ be a solution of \rf{state}. Then for any stopping times $0\les\t_1\les\t_2<\i$,
$$\ba{ll}
\ns\ds{1\over2}\dbE\lan P(\t_2,\a(\t_2))X(\t_2),X(\t_2)\ran-{1\over2}\dbE\lan P(\t_1,\a(\t_1)))X(\t_1),X(\t_1)\ran\\
\ns\ds={1\over2}\dbE\2n\int_{\t_1}^{\t_2}\3n\(\lan P_s(s,\a(s))X(s),X(s)\ran+\lan P(s,\a(s))X(s),A(\a(s))X(s)+B(\a(s))u(s)\ran\\
\ns\ds\q+{1\over2}\lan P(s,\a(s))[C(\a(s))X(s)+D(\a(s))u(s)],C(\a(s))X(s)\1n+\2n D(\a(s))u(s)\ran\\
\ns\ds\q+{1\over2}\L[\lan P(s,\a(s))X(s),X(s)\ran\)ds\\
\ns\ds={1\over2}\dbE\int_{\t_1}^{\t_2}\(\lan \[P_s(s,\a(s))+P(s,\a(s))A(\a(s))+A(\a(s))^\top P(s,\a(s))\\
\ns\ds\q+C(\a(s))^\top P(s,\a(s))C(\a(s))+\sum_{\jmath\in\cM} \l_{\imath\jmath}[P(s,\jmath)-P(s,\a(s))]\]X(s),X(s)\ran\\
\ns\ds\q+2\lan[D(\a(s))^\top P(s,\a(s))C(\a(s))+B(\a(s))^\top P(s,\a(s))]X(s),u(s)\ran\\
\ns\ds\q+\lan D(\a(s))^\top P(s,\a(s))D(\a(s))u(s),u(s)\ran\)ds.\ea$$

\ec

The proof is straightforward. Such a corollary will be useful later.

\section{Characterization of Cost Functional and Value Function}
\label{sec:cost}

\ms

In this section, we first look at representation of $J_{\scT}(t,x,\imath;u(\cd))$. To this end, we assume (A1) (and (A2)$'$ if $T=\i$). Since state equation \rf{state} is linear in $(X(\cd),u(\cd))$, one has
$$X(\cd)=X_0(\cd)+X_1(\cd),$$
with
$$\left\{\2n\ba{ll}
\ns\ds dX_0(s)=A(\a(s))X_0(s)ds+C(\a(s))X_0(s)dW(s),\qq s\in[t,T)],\\
\ns\ds X_0(t)=x,\ea\right.$$
and
$$\ba{ll}
\ns\ds\left\{\2n\ba{ll}
\ns\ds dX_1(s)=[A(\a(s))X_1(s)+B(\a(s))u(s)]ds\\
\ns\ds\qq\qq\qq+[C(\a(s))X_1(s)+D(\a(s))u(s)]
dW(s),\q s\in[t,T)],\\
\ns\ds X_1(t)=0.\ea\right.\ea$$
We may further write
$$X(\cd)=[\G_0x](\cd)+[\G_1u(\cd)](\cd),$$
for some linear and bounded operators $\G_0:\dbR^n\to\sX_{ad}[t,T]$ and $\G_1:\sU[t,T]\to\sX_{ad}[t,T]$. (in the case $0<T<\i$, (A2)$'$ is not necessary). The cost functional is purely quadratic terms. Thus, one has
$$\ba{ll}
\ns\ds J_{\scT}(t,x,\imath;u(\cd))={1\over2}\dbE\[\int_t^T\{\lan Q(\a(s))[X_0(s)+X_1(s)],X_0(s)+X_1(s)\ran\\
\ns\ds\qq\qq\qq\qq\qq+2\lan S(\a(s))[X_0(s)+X_1(s)],u(s)\ran+\lan R(\a(s))u(s),u(s)\ran\)ds\]\\
\ns\ds=\1n{1\over2}\dbE\[\2n\int_t^T\3n\(\lan Q(\a(s))X_1(s),X_1(s)\ran\1n+\1n2\lan S(\a(s))X_1(s),u(s)\ran\1n+\1n\lan R(\a(s))u(s),u(s)\ran\\
\ns\ds\qq\qq\q+2\lan Q(\a(s))X_0(s),X_1(t)\ran+2\lan S(\a(s))X_0(s),u(s)\ran\\
\ns\ds\qq\qq\q+\lan Q(\a(s))X_0(s),X_0(s)\ran\)ds\].\ea$$
Thus, we have
\bel{J=M_2x^2}J_{\scT}(t,x,\imath;u(\cd))=\lan M_{22}u,u\ran+2\lan M_{12}x,u\ran+\lan M_{11}x,x\ran,\ee
for some $M_{22}\in\cL(\sU[t,T])$ (the set of linear bounded operators from $\sU[t,T]$ to itself) with $M_{22}^*=M_{22}$, $M_{12}\in\cL(\dbR^n;\sU[t,T])$ (the set of linear bounded operators from $\dbR^n$ to $\sU[t,T]$), $M_{11}\in\dbR^{n\times n}$. This gives a representation of the cost functional. We now briefly look at the properties of the operators $M_{22}$, $M_{12}$ and $M_{11}$.

\ms

Note that (A1) (and (A2)$'$ if $T=\i$) only guarantees that for some $u(\cd)\in\sU[t,\i]$, $J(t,x,\imath;u(\cd))$ is well-defined. It quite be possible that
$$\inf_{u(\cd)\in\sU[t,T]}J(t,x,\imath;u(\cd))=-\i.$$
Here is a simple example.

\bex{} \rm Let $n=d=1$, $t<\i$, $\cM$ be a singleton, and consider controlled linear SDE
$$dX(s)=u(s)dW(s),\qq s\in[t,T],\qq X(t)=x.$$
The cost functional is
$$J_{\scT}(t,x;u(\cd))=\dbE\int_t^T\(|X(s)|^2-|u(s)|^2\)ds.$$
Then
$$\ba{ll}
\ns\ds J_{\scT}(t,x;u(\cd))=\dbE\int_t^T\(\big|x+\int_t^su(r)dW(r)\big|^2-|u(s)|^2\)ds\\
\ns\ds\qq\qq\q\;=|x|^2(T-t)+\dbE\int_t^T\[\(\int_t^s|u(r)|^2dr\)-|u(s)|^2\]
ds\\
\ns\ds\qq\qq\q\;=|x|^2(T-t)+\dbE\int_t^T(T-s-1)|u(s)|^2ds.\ea$$
Thus, under control $u_k(\cd)=k\chi_{[T-\d,T]}(\cd)$, with $\d\in(0,T-t)$, we have
$$J(t,x;u_k(\cd))=|x|^2(T-t)-{k\over2}\[1-(1-\d)^2\].$$
Hence,
$$\inf_{u(\cd)\in\sU[t,T]}J(t,x;u(\cd))=-\i.$$

\ex

The above example shows that unless we assume more, Problem (LQ)$_{\sc t,T}$ might be meaningless. In fact, for our LQ problem, we need to introduce the following definitions.

\bde{} \rm (i) Problem (LQ)$_{\sc t,T}$ is {\it finite} at $(t,x,\imath)\in\sD$ if
\bel{>-i}\inf_{u(\cd)\in\sU[t,T]}J_{\scT}(t,x,\imath;u(\cd))>-\i.\ee

\ms

(ii) Problem (LQ)$_{\sc t,T}$ is (uniquely) {\it open-loop solvable} at $(t,x,\imath)\in\sD$, if one can find a (the) control $\bar u^{t,x,\imath}_{\scT}(\cd)\in\sU[t,T]$ such that
\bel{J=J*}J_{\scT}(t,x,\imath;\bar u^{t,x,\imath}_{\scT}(\cd))=\inf_{u(\cd)\in\sU[t,T]}J_{\scT}(t,x,\imath;
u(\cd)).\ee

\ms

The following result is a basic result for quadratic functionals in Hilbert spaces.

\ede

\bp{} \sl Let \rf{J=M_2x^2} hold.

\ms

{\rm(i)} If \rf{>-i} holds, then
\bel{}M_{22}\ges0.\ee

\ms

{\rm(ii)} Function $u(\cd)\mapsto J(t,x,\imath:u(\cd))$ has a (unique) minimum if and only if 
\bel{M_12}M_{22}u+M_{12}x=0\ee
admits a (unique) solution $u(\cd)$. In particular, if
\bel{M_22>0}M_{22}\ges\d I,\ee
for some $\d>0$, then equation \rf{M_12} admits a unique solution
$$u=-M_{22}^{-1}M_{12}x,$$
and
\bel{V=}\ba{ll}
\ns\ds V_T(t,x,\imath)=\inf_{u(\cd)\in\sU[t,T]}J_T(t,x,\imath;u(\cd))=-\lan M_{22}^{-1}M_{12}x,M_{12}x\ran+\lan M_{11}x,x\ran\\
\ns\ds\qq\qq=\lan(M_{11}-M_{12}^*M_{22}^{-1}M_{12})x,x\ran={1\over2}\lan P(t,\imath)x,x\ran,\ea\ee
for some functions $P:[t,T)]\times\cM\to\dbS^n$.

\ep

It is easy to see that condition \rf{M_22>0} is equivalent to the uniform convexity of the map $u(\cd)\mapsto J_{\scT}(t,x,\imath;u(\cd))$, or for some $\d,K>0$, it holds
\bel{J>}J_{\scT}(t,x,\imath;u(\cd))\ges\d\int_t^T|u(s)|^2ds-K|x|^2,\qq\forall u(\cd)\in\sU_{ad}^{t,x,\imath}[t,T],\ee
Usually, condition \rf{M_22>0} is not easy to explicitly give. However, some sufficient conditions for that is easy to give. We now introduce the following assumption.

\ms

{\bf (A3)} Suppose that $Q(\imath)\in\dbS^n_{++}$, $R(\imath)\in\dbS^m_{++}$ such that
$$Q(\imath)-S(\imath)^\top R(\imath)^{-1}S(\imath)\in\dbS^n_{++}.$$

\ms

The following is mainly about the implication of (A3).

\bp{stconvex}\sl Let {\rm(A1)}, (also {\rm(A2)$'$} hold if $T=\i$), and {\rm(A3)} hold. Then the cost functional $u(\cd)\mapsto J_{\scT}(t,x,\imath;u(\cd))$ is uniformly convex over $\sU_{ad}^{t,x,\imath}[0,T]$, uniform in $0<T\les\i$.

\ep

\it Proof. \rm Note that for $0<\e<1$ small, and $a,b\in\dbR$, we have
$$|a|^2=|a+b-b|^2=|a+b|^2-2(a+b)b+|b|^2\les(1+\e)|a+b|^2+(1+{1\over\e})
|b|^2,$$
which leads to
$$|a+b|^2\ges{\frac1{1+\e}}|a|^2-{1\over\e}|b|^2.$$
Consequently, by suppressing $\imath\in\cM$,
$$\ba{ll}
\ns\ds2g(s,x,\imath,u)=\lan Qx,x\ran+2\lan Sx,u\ran+\lan Ru,u\ran\\
\ns\ds=\lan Qx,x\ran+2\lan R^{-{1\over2}}Sx,R^{1\over2}u\ran+|R^{1\over2}u|^2\\
\ns\ds=\lan[Q-S^\top R^{-1}S]x,x\ran+|R^{1\over2}u+R^{-{1\over2}}Sx|^2\\
\ns\ds\ges|R^{1\over2}u+R^{-{1\over2}}Sx|^2\ges{1\over1+\e}|R^{1\over2}
u|^2-{1\over\e}|R^{-{1\over}}Sx|^2\ea$$
Then we have \rf{J>} and the uniform convexity of the map $u(\cd)\mapsto J_{\scT}(t,x,\imath;u(\cd))$ on $\sU_{ad}^{t,x,\imath}[t,T]$ follows. It is clear that the convexity is uniform in $0<T\les\i$.  \endpf

\ms

On the other hand, under (A1) ((A2)$'$ if $T=\i$)), and (A3), we see that (for some small $\e>0$) we may write
$$\ba{ll}
\ns\ds Q(\imath)=(1-\e)Q(\imath)+\e Q(\imath),\q R(\imath)=(1-\e)R(\imath)+\e R(\imath),\\
\ns\ds(1-\e)Q(\imath)-{1\over1-\e}S(\imath)^\top R(\imath)^{-1}S(\imath)\in\dbS^n_{++},\ea$$
and
\bel{J^0>0}\3n\ba{ll}
\ns\ds2J_{\scT}(t,x,\imath;u(\cd))=\dbE\int_t^T2g(X(s),\a(s),u(s))ds\\
\ns\ds=\1n\dbE\2n\int_t^T\3n\(\big|[(1\1n-\1n\e)Q\1n-\1n{1\over1\1n-\1n\e}S^\top \1n R^{-1}\1n S]^{1\over2}X(s)\big|^2\2n+\1n\big|\sqrt{1\1n-\1n\e}R^{1\over2}u(s)\1n+\1n
{1\over\sqrt{1\1n-\1n\e}}
R^{-{1\over2}}SX(s)\big|^2\)ds\\
\ns\ds\qq\qq+\dbE\2n\int_t^T\3n\e\(\lan QX(s),X(s)\ran+\lan Ru(s)u(s)\ran\)ds\\
\ns\ds\ges\e\dbE\int_t^T\3n\(\lan QX(s),X(s)\ran+\lan Ru(s),u(s)\ran\)ds.\ea\ee
This implies that (see \rf{V=})
\bel{P>d}P(t,\imath)=2(M_{11}-M_{12}^*M_{22}^{-1}M_{12})\ges\d I,\ee
for some $\d>0$, independent of $t\ges0$ (only depends on (A3)).

\ms

We see that under (A1), (A2)$'$, $\sU^{0,x,\imath}[0,\i]=\sU[0,\i]$. Also, we note
\bel{QR}\ba{ll}
\ns\ds Q^{\Th^*}-(S^{\Th^*})^\top R^{-1}S^{\Th^*}\\
\ns\ds\qq\qq=Q+(\Th^*)^\top S+S^\top\Th^*+(\Th^*)^\top R\Th^*-(S+R\Th^*)^\top R^{-1}(S+R\Th^*)\\
\ns\ds\qq\qq=Q-S^\top R^{-1}S.\ea\ee
This means (A3) holds for $\{Q,S,R\}$ if and only if it holds for $\{Q^{\Th^*},S^{\Th^*},R\}$ for any $\Th^*((\cd)\in\hTh$.

\ms

Finally, from the well-known fact, we see that for the case $0<T<\i$, under (A1) and (A3), each Problem (LQ)$_{\sc t,\scT}$ is well-posed and admits a unique optimal control $\bar u_{\scT}(\cd)\in\sU[t,T]$; whereas under (A1), (A2)$'$, and (A3), Problem (LQ)$_{\sc0,\scT}$ admits a unique optimal control $\bar u_\i(\cd)\in\sU_{ad}^{0,x,\imath}[0,\i]$. We note that (A3) can be slightly relaxed, which leads to the so-called indefinite LQ problems. We omit the details here.

\ms

\section{ Main Results}\label{sec:hom}

Before going further, let us make another observation. We note that
\bel{J=J}\ba{ll}
\ns\ds J_{\sc\i}(t,x,\imath;u(\cd))=\dbE\int_t^\i g(X(s),\a(s),u(s))ds\\
\ns\ds=\dbE\int_0^\i g(X(s+t),\a(s+t),u(s+t))ds\\
\ns\ds=\dbE\int_0^\i g(\widehat X(s),\widehat\a(s),\widehat u(s))ds=\widehat J_{\scT}(0,x,\a(t);\widehat u(\cd)),\ea\ee
with
$$(\widehat X(s),\widehat\a(s),\widehat u(s))=(X(s+t),\a(s+t),u(s+t)),$$
which has the same distribution as that of $(X(s),\a(s),u(s))$ since the following holds:
$$\left\{\2n\ba{ll}
\ns\ds d\widehat X(s)=\(A(\widehat\a(s))\widehat X(s)+B(\widehat\a(s))\widehat u(s)\)ds+\(C(\widehat\a(s))\widehat X(s)+D(\widehat\a(s))\widehat u(s)\)d\widehat W(s),\\
\ns\ds\widehat X(0)=x,\q\widehat\a(0)=\imath,\ea\right.$$
where $s\mapsto\widehat W(s)\equiv W(s+t)$ is a Brownian motion starting from $W(t)$. then $(\widehat X(\cd),\widehat\a(\cd),\widehat u(\cd))$ and $(X(\cd),\a(\cd),u(\cd))$ might have different distributions, in general. (of course, assuming (A1) ((A2)$'$ if $T=\i$), and (A3)), we have the equivalence between Problems (LQ)$_{\sc t,\sc\i}$ and (LQ)$_{\sc0,\sc\i}$. Note that for Problem (LQ)$_{\sc0,\sc\i}$, the cost functional is actually equal to
(compare \rf{cost[t,T]} with $T=\i$)
\bel{J^0}J_{\sc\i}(x,\imath;u(\cd))=\dbE\int_0^\i g(X(s),\a(s),u(s))ds,\ee
which is independent of $t$. Consequently, the corresponding value function can be written as

\bel{V_T^0}V_{\scT}(x,\imath)=\inf_{u(\cd)\in\sU[0,\i]}J_{\scT}(x,\imath;u(\cd)).\ee

\ms

In order to study Problem (LQ)$_{\sc t,\scT}$ with $0\les t<T<\i$, let us make the following observation. First of all, from Corollary \ref{cor2.7}, we have
$$\ba{ll}
\ns\ds\dbE\lan P(T,\a(T))X(T),X(T)\ran-\dbE\lan P(t,\a(t))X(t),X(t)\ran\\
\ns\ds=\dbE\int_t^T\(\lan \[P_s(s,\a(s))+P(s,\a(s))A(\a(s))+A(\a(s))^\top P(s,\a(s))\\
\ns\ds\qq+C(\a(s))^\top P(s,\a(s))C(\a(s))+\L[P(s,\cd)](\a(s))\]X(s),X(s)\ran\\
\ns\ds\qq+2\lan[D(\a(s))^\top P(s,\a(s))C(\a(s))+B(\a(s))^\top P(s,\a(s))]X(s),u(s)\ran\\
\ns\ds\qq+\lan D(\a(s))^\top P(s,\a(s))D(\a(s))u(s),u(s)\ran\)ds.\ea$$
Thus, (suppressing $\a(\cd)$ and $s$)
$$\ba{ll}
\ns\ds J_{\scT}(t,x,\imath;u(\cd))={1\over2}\dbE\int_t^T\(\lan QX,X\ran+2\lan SX,u\ran+\lan Ru,u\ran\)ds\\
\ns\ds={1\over2}\lan P(t,\imath)x,x\ran+{\over2}\dbE\int_t^T\(\lan\big[\dot P+\L[P]+PA+A^\top P+C^\top PC+Q\big]X,X\ran\\
\ns\ds\qq\qq+2\lan(B^\top P+D^\top PC+S)X,u\ran+\lan(R+D^\top PD)u,u\ran \)ds\\
\ns\ds={1\over2}\lan P(t,\imath)x,x\ran+{1\over2}\dbE\int_t^T\(\lan\[\dot P+\L[P]+PA+A^\top P+C^\top PC+Q\\
\ns\ds\qq\qq-(B^\top\1n P\1n+\1n D^\top\1n PC\1n+\1n S)^\top\1n(R\1n+\1n D^\top\1n PD)^{-1}(B^\top\1n P\1n+\1n D^\top\1n PC\1n+\1n S)\]X,X\ran\\
\ns\ds\qq\qq+|(R+D^\top PD)^{1\over2}u+(R+D^\top PD)^{-{1\over2}}(B^\top P+D^\top PC+S)X|^2\)ds\\
\ns\ds\equiv{1\over2}\lan P(t,\imath)x,x\ran+{1\over2}\dbE\int_t^T\(\lan[\dot P+F(P)]X,X\ran+|(R+D^\top PD)^{1\over2}(u-\Th_TX)|^2\)ds,\ea$$
with an obvious definition of $F(P)$ and
\bel{Th_T}\Th_{\scT}(s,\imath)\1n=\1n-[R(\imath)\1n+\1n D(\imath)^\top\1n P_{\scT}(s,\imath)D(\imath)]^{-1}[B(\imath)^\top\1n P_{\scT}(s,\imath)\1n+\1n D(\imath)^\top\1n P_{\scT}(s,\imath)C(\imath)\1n+\1n S(\imath)].\ee
Hence, if we let $P_{\scT}:[t,T]\times\cM\to\dbS^n_{++}$ be the solution to the following {\it differential Riccati equation} (DRE, for short)  (suppressing $s\in[0,T]$ and $\imath\in\cM$ ):
\bel{DRE}\left\{\2n\ba{ll}
\ds\dot P_{\scT}+\L[P_{\scT}]+P_{\scT}A+A^\top P_{\scT}+C^\top P_{\scT}C+Q\\
\ns\ds-(P_{\scT}B\1n+\1n C^\top\1n P_{\scT}D\1n+\1n S^\top)(R\1n+\1n D^\top\1n P_{\scT}D)^{-1}(B^\top\1n P_{\scT}\1n+\1n D^\top\1n P_{\scT}C\1n+\1n S)\1n=\1n0,\q t\in[0,T],\\
\ns\ds P_{\scT}(T)=0,\ea\right.\ee
and control be the following state feedback form:
\bel{bar u}\bar u_{\scT}^{\sc t,x,\imath}(s)=\Th_{\scT}(s,\a(s))\bar X_{\scT}^{\sc t,x,\imath}(s),\qq s\in[t,T],\ee
then
\bel{V=PX^2}J_{\scT}(t,x,\imath;\bar u_{\scT}^{\sc t,x,\imath}(\cd))=V_{\scT}(t,x,\imath)={1\over2}\dbE\lan P_{\scT}(t;\imath)x,x\ran.\ee
This means the control given by \rf{bar u}, called a {\it closed-loop representation}, is optimal, and the corresponding state processes will be optimal as well. Hence, we obtain  the expressions of the optimal pair (of Problem (LQ)$_{\sc t,\scT}$). Thus, it is wise to introduce the differential Riccati equation \rf{DRE} when we study Problem (LQ)$_{\sc t,\scT}$. Now, a function $P_{\scT}(\cd\,,\cd)\in C([0,T],\dbS^n_{++})\times\cM$ is called a {\it $\d$-strong regular solution} of \rf{DRE} if it is a solution of \rf{DRE} and for the given $\d>0$,
\bel{R+DPD1}R(\imath)+D(\imath)^\top P_{\scT}(s,\imath)D(\imath)\ges\d I,\qq\forall(s,\imath)\in[0,T]\times\cM.\ee
From \cite{Zhang-2021}, we know that the above $\d$ is exactly the same constant appearing in the uniform convexity condition of the map $u(\cd)\mapsto J_{\scT}^{\sc0}(t,x,\imath;u(\cd))$. Also, with $\Th_{\scT}(s,\imath)$ given by \rf{Th_T}, we have that
\begin{align*}
& P_{\scT}A^{\Th_{\scT}}+(A^{\Th_{\scT}})^\top P_{\scT}+(C^{\Th_{\scT}})^\top P_{\scT}C^{\Th_{\scT}}+Q^{\Th_{\scT}}\\
&\qq-[P_{\scT}B+(C^{\Th_{\scT}})^\top P_{\scT}D+(S^{\Th_{\scT}})^\top](R+D^\top P_{\scT}D)^{-1}[B^\top P_{\scT}+D^\top P_{\scT}C^{\Th_{\scT}}+S^{\Th_{\scT}}]\\
&=P_{\scT}(A+B\Th_{\scT})+(A+B\Th_{\scT})^\top P_{\scT}+(C+D\Th_{\scT})^\top P_{\scT}(C+D\Th_{\scT})\\
&\qq+Q+(\Th_{\scT})^\top S+S^\top\Th_{\scT}+(\Th_{\scT})^\top R\Th_{\scT}\\
&\qq-[P_{\scT}B+(C+D\Th_{\scT})^\top P_{\scT}D+(S+R\Th_{\scT})^\top](R+D^\top P_{\scT}D)^{-1}\\
&\qq\qq\cd[B^\top P_{\scT}+D^\top P_{\scT}(C+D\Th_{\scT})+S+R\Th_{\scT}]\\
&=P_{\scT}A+P_{\scT}B\Th_{\scT}+A^\top P_{\scT}+(\Th_{\scT})^\top B^\top P_{\scT}+C^\top P_{\scT}C+(\Th_{\scT})^\top D^\top P_{\scT}D\Th_{\scT}\\
&\qq+C^\top P_{\scT}D\Th_{\scT}+(\Th_{\scT})^\top D^\top P_{\scT}C+Q+(\Th_{\scT})^\top S+S^\top\Th_{\scT}+(\Th_{\scT})^\top R\Th_{\scT}\\
&\qq-(P_{\scT}B+C^\top P_{\scT}D+S^\top)(R+D^\top P_{\scT}D)^{-1}(B^\top P_{\scT}+C^\top P_{\scT}D+S)\\
&\qq-(\Th_{\scT})^\top(D^\top P_{\scT}D+R)(R+D^\top P_{\scT}D)^{-1}(B^\top P_{\scT}+D^\top P_{\scT}C+S)\\
&\qq-(P_{\scT}B+C^\top P_{\scT}D+S)(R+D^\top P_{\scT}D)^{-1}(D^\top P_{\scT}D+R)\Th_{\scT}\\
&\qq-(\Th_{\scT})^\top(D^\top P_{\scT}D+R)(R+D^\top P_{\scT}D)^{-1}(D^\top P_{\scT}D+R)\Th_{\scT}\\
&=P_{\scT}A+P_{\scT}B\Th_{\scT}+A^\top P_{\scT}+(\Th_{\scT})^\top B^\top P_{\scT}+C^\top P_{\scT}C+Q\\
&\qq-(P_{\scT}B+C^\top P_{\scT}D+S^\top)(R+D^\top P_{\scT}D)^{-1}(B^\top P_{\scT}+C^\top P_{\scT}D+S)\\
&\qq+(\Th_{\scT})^\top D^\top P_{\scT}D\Th_{\scT}+C^\top P_{\scT}D\Th_{\scT}+(\Th_{\scT})^\top D^\top P_{\scT}C\\
&\qq-(\Th_{\scT})^\top(B^\top P_{\scT}+D^\top P_{\scT}C+S)-(P_{\scT}B+C^\top P_{\scT}D+S)\Th_{\scT}\\
&\qq-(\Th_{\scT})^\top(D^\top P_{\scT}D+R)\Th_{\scT}+(\Th_{\scT})^\top S+S^\top\Th_{\scT}+(\Th_{\scT})^\top R\Th_{\scT}\\
&=P_{\scT}A+A^\top P_{\scT}+C^\top P_{\scT}C+Q\\
&\qq-(P_{\scT}B+C^\top P_{\scT}D+S^\top)(R+D^\top P_{\scT}D)^{-1}(B^\top P_{\scT}+C^\top P_{\scT}D+S).
\end{align*}
Thus, differential Riccati equation \rf{DRE} can be written as
\bel{DRE*}\left\{\2n\ba{ll}
\ds\dot P_{\scT}+\L[P_{\scT}]+P_{\scT}A^{\Th_{\scT}}+(A^{\Th_{\scT}})^\top P_{\scT}+(C^{\Th_{\scT}})^\top P_{\scT}C^{\Th_{\scT}}+Q^{\Th_{\scT}}\\
\ns\ds-[P_{\scT}B^\top\2n+\1n(C^{\Th_{\scT}})^\top\1n P_{\scT}D+\1n (S^{\Th_{\scT}})\1n^\top\1n](R\1n+\1n D^\top\1n P_{\scT}D)^{-1}\1n(B^\top\2n P_{\scT}\1n+\1n D^\top\2n P_{\scT}C^{\Th_{\scT}}\2n+\2n S^{\Th_{\scT}})\1n=\1n0,\\
\ns\ds\qq\qq\qq\qq\qq\qq\qq\qq\qq\qq\qq\qq t\in[0,T],\\
\ns\ds P_{\scT}(T)=0,\ea\right.\ee
Likewise, for Problem (LQ)$_{\sc0,\sc\i}$, suggested by the above (and the formal limit of \rf{DRE}), we introduce the following {\it algebraic Riccati equation} (ARE, for short) for $P_{\sc\i}:\cM\to\dbS^n$ (suppressing $\imath\in\cM$ again):
\bel{ARE}\ba{ll}
\ns\ds\L[P_{\sc\i}]+P_{\sc\i}A+A^\top P{\sc\i}+C^\top P_{\sc\i}C+Q\\
\ns\ds\q-(B^\top P_{\sc\i}+D^\top P_{\sc\i}C+S)^\top(R+D^\top P_{\sc\i}D)^{-1}(B^\top P_{\sc\i}+D^\top P_{\sc\i}C+S)=0.\ea\ee
With the notation of $F(\cd)$, the ARE can be written as
\bel{ARE*}F(P_{\sc\i})=0.\ee
A map $P_{\sc\i}:\cM\to\dbS^n_{++}$ is called a {\it $\d$-strong regular solution} of \rf{ARE} if it is a solution of \rf{ARE}, and for the given $\d>0$,
\bel{R+DPD2}R(\imath)+D(\imath)^\top P_{\sc\i}(\imath)D(\imath)\ges\d I,\qq\forall\imath\in\cM.\ee
In addition, we require $\Th_{\sc\i}(\cd)\in\BS[A,C;B,D]$, with
\bel{DefThI}\Th_{\sc\i}(\imath)\1n=\1n-[R(\imath)\1n+\1n D(\imath)^\top P_{\sc\i}(\imath)D(\imath)]^{-1}[B(\imath)^\top P_{\sc\i}(\imath)\1n+\1n D(\imath)^\top P_{\sc\i}(\imath)C(\imath)\1n+\1n S(\imath)],\ee
The same as above, the ARE \rf{ARE} can be written as
\bel{ARE*-1}\ba{ll}
\ns\ds\L[P_{\sc\i}]+P_{\sc\i}A^{\Th_\i}+(A^{\Th_\i})^\top P{\sc\i}+(C^{\Th_\i})^\top P_{\sc\i}C^{\Th_\i}+Q^{\Th_\i}\\
\ns\ds-[P_{\sc\i}B+(C^{\Th_\i})^\top P_{\sc\i}D+(S^{\Th_\i})^\top](R+D^\top P_{\sc\i}D)^{-1}\\
\ns\ds\qq\cd[B^\top P_{\sc\i}+D^\top P_{\sc\i}C^{\Th_\i}+S^{\Th_\i}]=0.\ea\ee
We have seen that the given $\d>0$ measures certain level of the strong regularity of the solutions $P(\cd\,,\imath)$ and $P_{\sc\i}(\imath)$. The following lemma is about Problem (LQ)$_{\sc t,\scT}^{\sc0}$ with $0\les t<T\les\i$.

\bl{} \sl Let $0<T\les\i$. Let {\rm(A1)}, (also {\rm(A2)$'$} if $T=\i$), and {\rm(A3)} hold. Then

\ms

{\rm(i)} DRE \rf{DRE} admits a unique $\d$-strong regular solution $P_{\scT}(\cd\,,\cd)$ for some $\d>0$, uniform in $T>0$. Further,
\rf{V=PX^2} holds and the optimal control of Problem {\rm(LQ)$_{\sc t,\scT}$} admit the closed-loop representation \rf{bar u}.

\ms

{\rm(ii)} ARE \rf{ARE} admits a unique $\d$-strong regular solution $P_{\sc\i}(\cd)$ for some $\d>0$. Further,
\bel{V_i}V_{\sc\i}(x,\imath)={1\over2}\lan P_{\sc\i}(\imath)x,x\ran.\ee
and the optimal control of Problem {\rm(LQ)$_{\sc0,\sc\i}^{\sc0}$} admits the following closed-loop representation:
\bel{u_i}\bar u_{\sc\i}^{\sc0,x,\imath}(s)=\Th_{\sc\i}(\a(s))\bar X_{\sc\i}^{\sc0,x,\imath}(s),\qq s\in[0,\i),\ee
with $\Th_\i(\cd)$ given by \rf{DefThI}.

\ms

{\rm(iii)} For any given $t\in[0,\i)$, the follows holds
\bel{P to P}P_{\scT}(t,\imath)=P_{\scT-t}(0,\imath)\nearrow P_{\sc\i}(\imath),\qq\hb{as } T\nearrow\i,\q\forall\imath\in\cM.\ee
Moreover, there exists a $\d>0$ so that (for some $K,\d>0$)
\bel{expPi}
0\les P_{\sc\i}(\imath)-P_{\scT}(t,\imath)\les Ke^{-\d(T-t)}I,\qq t\in[0,T],\ee
and consequently,
\bel{expTh}|\Th_{\sc\i}(\imath)-\Th_{\scT}(t,\imath)|\les Ke^{-\d(T-t)},\qq t\in[0,T].\ee

\el

\it Proof. \rm (i). By Theorem 6.3 of \cite{Zhang-2021}, we have all the conclusions of (i).

\ms

(ii) and (iii). The proof is lengthy, and we split it into several steps.

\ms

\it Step 1. Equality in \rf{P to P}. its convergence, and the equation.

\ms

\rm Under (A1), (A2)$'$ and (A3), there exists a $P:\cM\to\dbS^n_{++}$ so that (See \rf{V=})
$$V_{\sc\i}(x,\imath)={1\over2}\lan P(\imath)x,x\ran.$$
Let $P(\cd\,,\cd)$ be the solution of \rf{DRE}. Define
$$\wt P(s,\imath)=P_{\scT}(s+t,\imath),\qq s\in[0,T-t],~\imath\in\cM,$$
Then we see that $\wt P(\cd\,,\cd)$ satisfies the same equation in \rf{DRE} with $\wt P(T-t)=0$. Thus, by the uniqueness, $\wt P(\cd\,,\cd)$ must be the same as $P_{\scT-t}(\cd\,,\cd)$ (They satisfy the same equation, and assume the same value at $T-t$).  Hence,
\bel{P=P}P_{\sc{T-t}}(s,\imath)=\wt P(s,\imath)=P_{\scT}(s+t,\imath),\qq s\in [0,T-t],~\imath\in\cM.\ee
In particular, (take $s=0$)
\bel{P=P2}P_{\sc{T-t}}(0,\imath)=P_{\scT}(t,\imath),\qq\imath\in\cM.\ee
Now, by (A3), it can be seen that $P_{\scT}(t,\imath)$ is increasing in $T$,
$$P_{\sc{T-t}}(0,\imath)=P_{\scT}(t,\imath)\les P(\imath).$$
Therefore there exists a $\wt P(\imath)$ independent of $t$ such that
\bel{P to P*}\wt P(\imath)=\lim_{T\to\i}P_{\scT}(t,\imath).\ee
From \rf{DRE}, for any given $0\les\t_1<\t_2<T\les\i$,
$$P_{\scT}(\t_2,\imath)-P_{\scT}(\t_1,\imath)=\int_{\t_1}^{\t_2}F(P_{\scT}(s,
\imath),\a(\imath))ds,\qq \imath\in\cM,$$
for some $F(\cd\,,\cd)$. Taking $T\to\i$ in the above, we see that
$$0=(\t_2-\t_1)F(\wt P(\imath),\a(\imath)),\qq\imath\in\cM.$$
Thus, $\wt P(\cd)$ is a solution to \rf{ARE}.

\ms

Define
$$\wt\Th=-[R+D^\top\wt PD]^{-1}[B^\top\wt P+D^\top\wt PC+S],$$
with $\imath\in\cM$ suppressed.

\ms

\it Step 2. $\wt\Th\in\BS[A,C;B,D]$.

\ms

\rm
We have for any $T_0>0$,
\bel{convergence}\lim_{T\to\i}\sup_{0\les s\les T_0}|\Th_{\scT}(s,\imath)-\wt\Th(\imath)|=0,\ee
where $\Th_{\scT}(\cd\,,\cd)$ is given by \rf{Th_T}.
Now let us show $\wt\Th(\cd)\in\BS[A,C;B,D]$. To this end, let us observe the following system:
\bel{state6}\left\{\2n\ba{ll}
\ds d\wt X_{\sc\i}(s)=[A(\a(s))+B(\a(s))\wt\Th(\a(s))]\wt X_{\sc\i}(s)ds\\
\ns\ds\qq\qq\qq+[C(\a(s))+D(\a(s))\wt\Th(\a(s))\wt X_{\sc\i}(s)dW(s),\q s\in[0,\i),\\
\ns\ds\wt X_{\sc\i}(0)=x,\q\a(0)=\imath,\ea\right.\ee
This system is the homogeneous system under state feedback control
\bel{wtu}\wt u_{\sc\i}(s)=\wt\Th(\a(s))\wt X_{\sc\i}(s),\qq s\in[0,\i).\ee
Now, by the convergence \rf{convergence}, for any given $T_0>0$, we have
$$\lim_{T\to\i}\dbE\[\sup_{s\in[0,T_0]}|\bar X_{\scT}^{t,x,\imath}(s)-\wt X_{\sc\i}(s)|^2\]=0,$$
which leads to
$$\lim_{T\to0}\dbE\[\sup_{s\in[0,T_0]}|\bar u^{t,x,\imath}_{\scT}(s)-\wt u_{\sc\i}(s)|^2\]=0.$$
Note that under (A6), from \rf{J^0>0}, one has. for some $\e>0$
$$J_{\scT}(0,x,\imath;u(\cd))=\e\dbE\int_0^T\(|X(s)|^2+|u(s)|^2\)
ds,$$
for any state-control pair $(X(\cd),u(\cd))$. In particular,
$$\dbE\2n\int_0^{T_0}\3n\(|\bar X^{t,x,\imath}_{\scT}(s)|^2\2n+\1n|\bar u^{t,x,\imath}_{\scT}(s)|^2\)ds\1n\les\1n{1\over\e}J_{\scT}(t,x,\imath;
\bar u^{t,x,\imath}_{\scT}(\cd))\1n=\1n{1\over2\e}\lan P(t,\imath)x,x\ran\1n\les\1n{1\over2\e}\lan P(\imath)x,x\ran.$$
Thus, by Fatou's Lemma, we have
$$\int_0^{T_0}\(|\wt X_{\sc\i}(s)|^2+|\wt u_{\sc\i}(s)|^2\)ds\les{1\over2\e}\lan P(\imath)x,x\ran.$$
Since $T_0>0$ is arbitrary, we have $\wt u_{\sc\i}(\cd)\in\sU[0,\i]$, and $\wt X_{\sc\i}(\cd)\in\sX[0,\i]$. Thus, $\wt \Th(\cd)\in\BS[A,C;B,D]$.

\ms

\it Step 3. $\wt P(\imath)=P_{\sc\i}(\imath)$.

\ms

\rm

Let $X(\cd)$ be the solution to \eqref{state3} under $u(\cd)$. Then applying the extended It\^o's formula to $s\mapsto\lan\wt P(\a(s)) X(s),X(s)\ran$ (from Corollary 2.6 and (A6)), one can conclude that
$$\ba{ll}
\ns\ds d\lan\wt P(\a(s))X(s),X(s)\ran=\(\lan (\L[\wt P]+A^\top\wt P+\wt PA+C^\top\wt PC) X,X\ran\\
\ns\ds\qq\qq\qq\q+\lan D^\top\wt PD u,u\ran+2\lan(B^\top\wt P +D^\top\wt PC)X,u\ran\)ds+\{\cds\}dW(s).\ea$$
Thus, similar to the procedure right before \rf{DRE}, we have
$$J_{\sc\i}(0,x,\imath;u(\cd))=\lan \wt P(\imath)x,x\ran+\dbE\int_0^\i\lan(R+D^\top\wt PD)(u-\wt\Th X),u-\wt\Th X\ran ds.$$
Because $\wt\Th(\cd)\in\BS[A,C;B,D]$, we see from the above that \rf{wtu} gives the optimal control (of state feedback form). Hence, it is necessary that
$$\wt\Th(\a(s))=\Th_{\sc\i}(\a(s)),\qq s\in[0,\i).$$
Consequently,
$$\wt P(\imath)=P_{\sc\i}(\imath),\qq\imath\in\cM.$$
Then \rf{P to P} follows from \rf{P to P*}.

\ms

\it Step 4. Proof of \rf{expPi}--\rf{expTh}.

\ms

\rm

We first let $(\bar X^{\sc t,x,\imath}_{\scT}(\cd),\bar u^{\sc t,x,\imath}_{\scT}(\cd))$ be an optimal pair of Problem (LQ)$_{\sc t,\scT}^{\sc0}$. Then
$$\bar u_{\scT}^{\sc t,x,\imath}(s)=\Th_{\scT}(s,
\a(s))\bar X_{\scT}^{t,x,\imath}(s).$$
Apply the It\^o's formula to the map $s\mapsto\lan P_{\scT}(s,\a(s))\bar X_{\scT}^{\sc t,x,\imath}(s),\bar X_{\scT}^{\sc t,x,\imath}(s)\ran$ on $[t,T]$,
suppressing $\a(\cd)$ and $(t,x,\imath)$.
$$\ba{ll}
\ns\ds d\Big\{\dbE\lan P_{\scT}(s)\bar X_{\scT}(s),\bar X_{\scT}(s)\ran\Big\}\\
\ns\ds=\dbE\Big\{\lan\big(\dot P_{\scT}+\L[P_{\scT}]+P_{\scT}A+A^\top P_{\scT}+C^\top P_{\scT}C\big)\bar X_{\scT}(s),\bar X_{\scT}(s)\ran\\
\ns\ds\qq+2\lan\big(B^\top P_{\scT}+D^\top P_{\scT}C\big)\bar X_{\scT}(s),\bar u_{\scT}(s)\ran+\lan D^\top P_{\scT}D\bar u_{\scT}(s),\bar u_{\scT}(s)\ran\Big\}ds\\
\ns\ds=\dbE\Big\{\lan\big(\dot P_{\scT}+\L[P_{\scT}]+P_{\scT}A+A^\top P_{\scT}+C^\top P_{\scT}C+Q\big)\bar X_{\scT}(s),\bar X_{\scT}(s)\ran\\
\ns\ds\qq+2\lan\big(B^\top\1n P_{\scT}+D^\top\1n P_{\scT}C\1n+\1n S^\top\1n\big)\bar X_{\scT}(s),\bar u_{\scT}(s)\ran\1n+\1n\lan [R\1n+\1n D^\top\1n P_{\scT}D]\bar u_{\scT}(s),\bar u_{\scT}(s)\ran\\
\ns\ds\qq-\lan Q\bar X_{\scT}(s),\bar X_{\scT}(s)\ran-2\lan S^\top\bar X_{\scT}(s),\bar u_{\scT}(s)\ran-\lan R\bar u_{\scT}(s),\bar u_{\scT}(s)\ran\Big\}ds\\
\ns\ds=\dbE\Big\{\lan\big(\dot P_{\scT}(s)+\L[P_{\scT}(s)]+P_{\scT}(s)A+A^\top P_{\scT}(s)+C^\top P_{\scT}(s)C+Q\\
\ns\ds\qq-(P_{\scT}B\1n+\1n C^\top\1n P_{\scT}D\1n+\1n S^\top)(R\1n+\1n D^\top\1n P_{\scT}D)^{-1}(B^\top\1n P_{\scT}\1n+\1n D^\top P_{\scT}C\1n+\1n S)\big)\bar X_{\scT}(s),\bar X_{\scT}(s)\ran\\
\ns\ds\qq+|(R+D^\top P_{\scT}D)^{1\over2}\bar u_{\scT}(s)+(R+D^\top P_{\scT}D)^{-{1\over2}}(B^\top P_{\scT}+D^\top P_{\scT}C+S)\bar X_{\scT}(s)|^2\\
\ns\ds\qq-\lan Q\bar X_{\scT}(s),\bar X_{\scT}(s)\ran-2\lan S^\top\bar X_{\scT}(s),\bar u_{\scT}(s)\ran-\lan R\bar u_{\scT}(s),\bar u_{\scT}(s)\ran\Big\}ds\\
\ns\ds\les-\e\dbE|\bar X_{\scT}(s)|^2\les-\d\dbE\lan P_{\scT}(s)\bar X_{\scT}(s),\bar X_{\scT}(s)\ran.\ea$$
By Gronwall's inequality,
$$\dbE\lan P_{\scT}(s,\a(s))\bar X_{\scT}(s),\bar X_{\scT}(s)\ran\les Ke^{-\d(s-t)}\lan P_{\scT}(t,\imath)x,x\ran,\qq s\in[t,T].$$
Hence, with possibly different constants $K,\d>0$, one has
\bel{EX}\dbE|\bar X_{\scT}(s)|^2\les K|x|^2e^{-\d(s-t)},\qq s\in[t,T].\ee
By the dynamic program principle, we have
$$\ba{ll}
\ns\ds{1\over2}\lan P_{\sc\i}(\imath)x,x\ran=V_{\sc\i}(x,\imath)=\inf_{u(\cd)\in\sU[0,\i]}
J_{\sc\i}(0,x,\imath;u(\cd))=\inf_{u(\cd)\in\sU[t,\i]}
J_{\sc\i}(t,x,\imath;u(\cd))\\
\ns\ds=\inf_{u(\cd)\in\sU[t,T]}\dbE\(\int_t^Tg(X(s),\a(s),u(s))ds
+V_{\scT}(X(T),\a(T))\)\\
\ns\ds=\inf_{u(\cd)\in\sU[0,T]}\Big[J_{\scT }(0,x,\imath;u(\cd))+\dbE\lan P_{\scT}(T,\a(T)) X(T),X(T)\ran\Big].\ea$$
Consequently,
$$\ba{ll}
\ns\ds{1\over2}\lan P_{\scT}(t,\imath) x,x\ran=\inf_{u(\cd)\in \sU[t,T]}J_{\scT}(t,x,\imath;u(\cd))=J_{\scT}(t,\imath;
\bar u_{\scT}^{\sc t,x,\imath}(\cd))\\
\ns\ds=\Big[J_{\scT}(t,x,\imath;\bar u_{\scT}^{\sc t,x,\imath}(\cd))+\dbE\lan P_{\scT}(T,\a(T))\bar X_{\scT}^{\sc t,x,\imath}(T),\bar X_{\scT}^{\sc t,x,\imath}(T))\ran\Big]\\
\ns\ds\q-\dbE\lan P_{\scT}(T,\a(T))\bar X_{\scT}^{\sc t,x,\imath}(T),\bar X_{\scT}^{\sc t,x,\imath}(T)\ran\\
\ns\ds\ges\inf_{u(\cd)\in\sU[0,T]}\Big[J_{\scT }(0,x,\imath;u(\cd))+\dbE\lan P_{\scT}(T,\a(T)) X(T),X(T)\ran\Big]-K|x|^2e^{-\d(T-t)}\\
\ns\ds=\inf_{u(\cd)\in \sU[0,\infty]}J_{\sc\i}(0,x,\imath;u(\cd))-K|x|^2e^{-\d(T-t)}
=V_{\sc\i}
(x,\imath)-K|x|^2e^{-\d(T-t)}\\
\ns\ds={1\over2}\lan P_{\sc\i}(
\imath)x,x\ran-K|x|^2e^{-\d(T-t)}.\ea$$
Hence,
$$0\les\lan [P_{\sc\i}(\imath)-P_{\scT}(t,\imath)]x,x\ran\les K|x|^2e^{-\d(T-t)},$$
proving \rf{expPi}. Finally, by the representation of $\Th_\i(\cd)$ in \rf{DefThI} and $\Th_{\scT}(\cd)$ in \rf{Th_T}, and noting the fact that $P_{\scT}(\cd\,,\cd)$ and $P_{\sc\i}(\cd)$ are the $\d$-strong regular solutions of \rf{DRE} and \rf{ARE}, respectively, we have
$$\ba{ll}
\ns\ds|\Th_{\sc\i}(\imath)-\Th_{\scT}(t,\imath)|\les{1\over\d}\big|B(\imath)^\top
[P_{\sc\i}(\imath)-P_{\scT}(t,\imath)]+D(\imath)^\top[P_{\sc\i}(\imath)-
P_{\scT}(t,\imath)]C(\imath)\big|\\
\ns\ds\qq+{1\over\d^2}\big|D(\imath)^\top[P_{\sc\i}(\imath)-P_{\scT}(t,\imath)]
D(\imath)\big|\big|B(\imath)^\top P_{\sc\i}(\imath)+D(\imath)^\top P_{\scT}(\imath)C(\imath)+S(\imath)\big|\\
\ns\ds\qq\les Ke^{-\d(T-t)}.\ea$$
This proves \rf{expTh}. \endpf

\ms

The above proof of \rf{expPi} is different from that of Theorem 4.4 in \cite{Sun-Wang-Yong-2022}, which plays an essential role in studying STP for LQ optimal control problem without switching.
Now we are ready to present the main result of this paper.

\bt{} \sl Let {\rm(A1), (A2)$'$, (A3)} hold. Then there exists $\d>0$ such that for any $(t,x_{\scT},\imath),(t,x_{\sc\i},\imath)\in\sD$,
\bel{Turncase0eq}\ba{ll}
\ns\ds\dbE\Big(|\bar X_{\scT}^{t,\imath,x_{\scT}}(s)-\bar X_{\sc\i}^{t,\imath,x_\i}(s)|^2+|\bar u_{\scT}^{t,\imath,x_T}(s)-\bar u_{\sc\i}^{t,\imath,x_\i}(s)|^2\Big)\\
\ns\ds\les Ke^{-\d(s-t)}|x_{\scT}-x_{\sc\i}|^2+Ke^{-\d(s-t)}
e^{-2\d(T-s)}|x_{\scT}|^2,\qq s\in[t,T].\ea\ee
Moreover, for $x_{\scT}=x_{\sc\i}=x$, $t=0$,
\bel{Eint}\limsup_{T\to\i}\dbE\int_0^T \Big(|\bar X_{\scT}^{0,\imath,x}(s)-\bar X_{\sc\i}^{0,\imath,x}(s)|^2+|\bar u_{\scT}^{0,\imath,x}(s)-\bar u_{\sc\i}^{0,\imath,x}(s)|^2\Big)ds=0.\ee

\et

\it Proof. \rm Let us suppress the superscripts $(t,x_{\scT},\imath)$ and $(t,x_{\sc\i},\imath)$ in the following. For reader's convenience, we rewrite two systems:
$$\left\{\2n\ba{ll}
\ns\ds d\bar X_{\scT}(s)=A^{\Th_T}\bar X_{\scT}(s)ds+C^{\Th_T}\bar X_{\scT}(s)dW(s),\qq s\in[t,T].\\
\ns\ds\bar X_{\scT}(t)=x_{\scT},\q\a(t)=\imath,\ea\right.$$
$$\left\{\2n\ba{ll}
\ns\ds d\bar X_{\sc\i}(s)=A^{\Th_\i}\bar X_{\sc\i}(s)ds+C^{\Th_\i}\bar X_{\sc\i}(s)dW(s),\qq s\in[t,\i),\\
\ns\ds\bar X_{\sc\i}(t)=x_{\sc\i},\q\a(t)=\imath.\ea\right.$$
Then
$$\left\{\2n\ba{ll}
\ns\ds d[\bar X_{\scT}(s)-\bar X_{\sc\i}(s)]=[A^{\Th_T}\bar X_{\scT}(s)-A^{\Th_\i}\bar X_{\sc\i}(s)]ds\\
\ns\ds\qq\qq\qq\qq\qq+[C^{\Th_T}\bar X_{\scT}(s)-C^{\Th_\i}\bar X_{\sc\i}(s)]dW(s),\qq s\in[t,\i),\\
\ns\ds\bar X_{\scT}(t)-\bar X_{\sc\i}(t)=x_{\scT}-x_{\sc\i},\q\a(t)=\imath.\ea\right.$$
Note that system $[A^{\Th_\i},C^{\Th_\i}]$ is stable and hence is dissipative. Therefore, for some $\Si(\cd)\in\hTh$ and $\bar\d>0$,
$$\(\L[\Si]+\Si A^{\Th_\i}+(A^{\Th_\i})^\top\Si+(C^{\Th_\i})^\top\Si C^{\Th_\i}\)(\imath)\les-2\d\Si(\imath),\q\imath\in\cM.$$
The above can be achieved if we shrink $\d>0$. Applying It\^o's formula to the map
$$s\mapsto\dbE\lan\Si(\a(s))[\bar X_{\scT}(s)-\bar X_{\sc\i}(s)],\bar X_{\scT}(s)-\bar X_{\sc\i}(s)\ran,$$
one has
$$\ba{ll}
\ns\ds{d\over ds}\dbE\lan\Si(\a(s))[\bar X_{\scT}(s)-\bar X_{\sc\i}(s)],\bar X_{\scT}(s)-\bar X_{\sc\i}(s)\ran\\
\ns\ds=\dbE\[\lan\L[\Si](\bar X_{\scT}-\bar X_{\sc\i}),\bar X_{\scT}-\bar X_{\sc\i}\ran+\lan\Si(A^{\Th_{\scT}}\bar X_{\scT}(s)-A^{\Th_\i}\bar X_{\sc\i},\bar X_{\scT}-\bar X_{\sc\i}\ran\\
\ns\ds\qq+\lan\Si(\bar X_{\scT}-\bar X_{\sc\i}),A^{\Th_{\scT}}\bar X_{\scT}-A^{\Th_\i}\bar X_{\sc\i}\ran\\
\ns\ds\qq+\lan\Si(C^{\Th_{\scT}}\bar X_{\scT}-C^{\Th_\i}\bar X_{\sc\i}),C^{\Th_{\scT}}\bar X_{\scT}-C^{\Th_\i}\bar X_{\sc\i}\ran\]\\
\ns\ds=\dbE\[\lan\(\L[\Si]+\Si A^{\Th_\i}+(A^{\Th_\i})^\top\Si+(C^{\Th_\i})^\top\Si C^{\Th_\i}\)(\bar X_{\scT}-\bar X_{\sc\i}),\bar X_{\scT}-\bar X_{\sc\i}\Big\ran\\
\ns\ds\qq+\lan\Si(A^{\Th_{\scT}}-A^{\Th_\i})\bar X_{\scT},\bar X_{\scT}-\bar X_\i\ran+\lan\Si(\bar X_{\scT}-\bar X_{\sc\i}),(A^{\Th_{\scT}}-A^{\Th_\i})\bar X_{\scT}\ran\\
\ns\ds\qq+\lan\Si(C^{\Th_{\scT}}-C^{\Th_\i})\bar X_{\scT},(C^{\Th_{\scT}}-C^{\Th_\i})\bar X_{\scT}\ran\\
\ns\ds\qq+2\lan\Si C^{\Th_\i}(\bar X_{\scT}-\bar X_{\scT}),
(C^{\Th_T}-C^{\Th_\i})\bar X_{\scT}\ran\]\\
\ns\ds\les-2\d\dbE|\Si^{1\over2}(\bar X_{\scT}-\bar X_{\sc\i})|^2\\
\ns\ds\qq+K|\Th_{\scT}-\Th_{\sc\i}|\dbE\(|\Si^{1\over2}(\bar X_{\scT}-\bar X_{\sc\i})||\bar X_{\scT}|\)+K|\Th_{\scT}-\Th_{\sc\i}|^2\dbE|\bar X_{\scT}|^2\\
\ns\ds\les-\d\dbE|\Si^{1\over2}(\bar X_{\scT}-\bar X_{\sc\i})|^2+Ke^{-2\d(T-s)}|x_{\scT}|^2
e^{-\d(s-t)},\\
\ns\ds\les-\d\dbE|\Si^{1\over2}(\bar X_{\scT}-\bar X_{\sc\i})|^2+K|x_{\scT}|^2
e^{-2\d(T-s)}e^{-\d(s-t)}.\ea$$
By Grownwall's inequality, we have
$$\ba{ll}
\ns\ds\dbE|\Si^{1\over2}(\a(s))[\bar X_{\scT}(s)-\bar X_{\sc\i}(s)]|^2\\
\ns\ds\les e^{-\d(s-t)}\dbE|\Si^{1\over2}(\imath)(x_{\scT}-x_{\sc\i})|^2
+K|x_{\scT}|^2\int_t^se^{-\d(s-r)}e^{-2\d(T-r)}e^{-\d(r-t)}dr.\ea$$
Thus,
$$\dbE|\bar X_{\scT}(s)-\bar X_{\sc\i}(s)|^2\les Ke^{-\d(s-t)}|x_{\scT}-x_{\sc\i}|^2+Ke^{-\d(s-t)}
e^{-2\d(T-s)}|x_{\scT}|^2.$$
Now, by \rf{bar u} and \rf{u_i}, we have
$$\ba{ll}
\ns\ds\dbE|\bar u_{\scT}(s)-\bar u_{\sc\i}(s)|^2=\dbE|\Th_{\scT}\bar X_{\scT}(s)-\Th_{\sc\i}\bar X_{\sc\i}(s)|^2\\
\ns\ds\les2|\Th_{\scT}-\Th_{\sc\i}|^2\dbE|\bar X_{\scT}|^2
+2|\Th_{\sc\i}|^2\dbE|\bar X_{\scT}(s)-\bar X_{\sc\i}(s)|^2\\
\ns\ds\les Ke^{-\d(T-s)}e^{-\d(s-t)}|x_{\scT}|^2+K\[Ke^{-\d(s-t)}|x_{\scT}-x_{\sc\i}|^2+Ke^{-\d(s-t)}
e^{-2\d(T-s)}|x_{\scT}|^2\].\ea$$
Thus, \rf{Turncase0eq} follows. Finally, \rf{Eint} follows very straightforwardly. \endpf


\section{Concluding Remarks}\label{sec:con}

In this paper, we have established the turnpike property for LQ optimal control in an infinite horizon with a regime-switching state, focusing on homogeneous cases. The results obtained will play an essential role in future research on non-homogeneous cases.   The method proposed to tackle the convergence of AREs provides some new insights for future works on turnpike properties on more complicated LQ optimal control problems.

\end{document}